\documentclass[11pt]{amsart}

\usepackage{graphicx,color,latexsym}
\usepackage{amssymb}
\usepackage{pb-diagram}
\usepackage{verbatim}

\newtheorem{theorem}{Theorem}[section]

\newtheorem{lemma}[theorem]{Lemma}
\newtheorem{remark}[theorem]{Remark}

\newtheorem{definition}[theorem]{Definition}
\newtheorem{cor}[theorem]{Corollary}
\newtheorem{conjecture}[theorem]{Conjecture}

\newtheorem{example}[theorem]{Example}

\begin{document}

\title{The Relative Symplectic Cone and $T^2-$fibrations}

\author{Josef G. Dorfmeister}
\address{School  of Mathematics\\  University of Minnesota\\ Minneapolis, MN 55455}
\email{dorfmeis@math.umn.edu}

\author{Tian-Jun Li}
\address{School  of Mathematics\\  University of Minnesota\\ Minneapolis, MN 55455}
\email{tjli@math.umn.edu}

\date{\today}

\begin{abstract}
In this note we introduce the notion of the relative symplectic cone.  As an application, we determine the {\it symplectic cone} of certain $T^2$-fibrations.  In particular, for some elliptic surfaces we verify a conjecture on the symplectic cone of minimal K\"ahler surfaces raised in \cite{TJL1}.

\end{abstract}

\maketitle

\tableofcontents

\section{Introduction}
Given an oriented smooth manifold $M$ known to admit symplectic
structures, one would like to know which cohomology classes $\alpha\in H^2(M,\mathbb R)$ can be
represented by an orientation compatible symplectic form $\omega\in\Omega^2(M)$ .  We shall always restrict ourselves to symplectic forms which are compatible with the fixed orientation of the manifold $M$.  This leads naturally to the
definition of the symplectic cone:
\begin{equation}
{\mathcal C}_M=\{\alpha\in H^2(M)\;\vert\;[\omega]=\alpha,\;  \omega \mbox{
is a symplectic form on }
M\}.
\end{equation}
In dimension 4, the symplectic cone has been determined in the
following cases: $S^2$-bundles (\cite{M2}), $T^2$-bundles over $T^2$
 (\cite{Ge}), all $b^+=1$ manifolds
(\cite{TJLL}), and minimal manifolds underlying a K\"ahler surface with
Kodaira dimension 0 (\cite{TJL1}).  A smooth  $4$-manifold $M$
is said to be minimal if it contains no exceptional class, i.e. a degree 2
homology class represented by a smoothly embedded sphere of self
intersection $-1$.

 Clearly, $ {\mathcal C}_M$ is contained in ${\mathcal P}_M$, the cone of
classes of positive squares in $H^2(M)$. Amazingly, when $M$ is a minimal representative of one of the previous examples, the symplectic cone
${\mathcal C}_M$ is actually equal to ${\mathcal P}_M$.

In general, ${\mathcal C}_M$ is smaller than ${\mathcal P}_M$, as there are
constraints coming from the Seiberg-Witten basic classes.  This is a consequence of
Taubes' remarkable equivalence between Seiberg Witten invariants SW and  Gromov invariants Gr (\cite{T}). As
exceptional classes and the canonical class of any symplectic structure all
give rise to SW basic classes, there are corresponding constraints
on $\mathcal C_M$.

If the smooth manifold $M$ underlies a minimal K\"ahler surface, a basic fact
(\cite{Witten}, \cite{FM}) is that all symplectic structures on $M$
have the same canonical class up to sign. Denote and fix one such
choice $-c_1(M)$.
 In light of this beautiful fact,  in \cite{TJL1}, the following conjecture was raised:

\begin{conjecture} \label{conj} If $M$ underlies a {\bf minimal}
K\"ahler surface, the symplectic cone ${\mathcal C}_M$ is equal to
${\mathcal P}^{c_1(M)}\cup{\mathcal P}^{-c_1(M)}$.\end{conjecture}

We define ${\mathcal P}^{\alpha}=\{e\in {\mathcal P}_M|e\cdot \alpha> 0 \}$
for nonzero $\alpha\in H^2(M;{\mathbb R})$ and ${\mathcal
P}^0={\mathcal P}_M$. As ${\mathcal P}^0={\mathcal P}_M={\mathcal
C}_M$ when $b^+(M)=1$ or $c_1(M)=0$ (i.e. $c_1(M)$ is torsion), this conjecture is known to be
true when $M$ underlies a minimal K\"ahler surface with $p_g=0$ or
Kodaira dimension $0$ (see also \cite{LU} for the case $p_g>0$).  Moreover, for any minimal K\"ahler surface, ${\mathcal C}_M\subset {\mathcal P}^{c_1(M)}\cup{\mathcal P}^{-c_1(M)}$.

In this note we will  show that this conjecture holds for certain
manifolds underlying minimal K\"ahler manifolds with $p_g>0$ and
Kodaira dimension $\kappa= 1$. Many such manifolds are
$T^2-$fibrations and can be written as a $T^2-$fiber sum of
manifolds with $p_g=0$ or Kodaira dimension 0.

 There are many ways to explicitly construct new
symplectic manifolds.  Common among most of these methods is that some type of surgery is
performed with respect to a codimension 2 symplectic submanifold $V$. It
is natural to ask how the symplectic forms on the new manifolds
relate to those on the constituent manifolds. This leads naturally
to the notion of the relative symplectic cone ${\mathcal C}_M^V$ defined
in Section \ref{relcone}.  As examples, we consider $T^2$ fibrations over $T^2$ (see \cite{Ge}) and manifolds with $b^+=1$.  These are of interest, as we will consider $T^2$ fiber sums in the following sections.

 The fiber sum of symplectic manifolds $X$ and $Y$
along symplectic embeddings of a codimension 2  symplectic manifold
$V$, denoted $M=X\#_VY$, as defined by Gompf (\cite{Go2}) and
McCarthy-Wolfson (\cite{MW}), and its inverse operation, the
symplectic cut, defined by Lerman (\cite{Le}), are briefly described
in Section \ref{cutsum}.

 We then proceed to show that  the sum and
cut operations naturally describe a cone $\mathcal C^{sum}$ of sum forms
in terms of the relative cones of $X$ and $Y$ with respect to $V$.  We also observe that in the case $V$ having trivial normal bundle,
${\mathcal C}^{sum}$ is actually  a subcone of the relative cone ${\mathcal C}^V_M$.

Furthermore, under some  topological restrictions on the sum
$M=X\#_VY$ and the respective relative cones ${\mathcal C}_*^V$, we show that the {\it relative}
symplectic cone ${\mathcal C}^V_M$ is actually equal to this subcone.

What does this imply for the symplectic cone of $M$?
 Notice that for a minimal
$T^2$-fibration, the canonical class is proportional to the fiber
class. Thus the relative symplectic cone, which is of course contained in the symplectic cone,  is essentially equal to the
symplectic cone.  This strategy applies perfectly to  fiber sums where one summand is
a product $T^2-$fibration, hence verifying the conjecture for such
$T^2$ fibrations.  During the preparation of this paper,  Friedl and Vidussi (see \cite{FV1} and \cite{FV2}) obtained a result, which allows them to determine the symplectic cone of $T^2\times \Sigma_g$.  They consider manifolds admitting a free circle action and use properties of the  Thurston norm ball of a quotient 3-manifold.  Their results allow them to determine completely the symplectic cone of a 4-manifold with a free circle action such that the orbit space is a graph manifold.   

We include an appendix concerning genericity results for almost complex structures $J$ which make $V$ pseudoholomorphic.  These results are needed to determine the relative cone in the $b^+=1$ case of Section \ref{relcone}.  They show that the set $\mathcal J_V$ of such almost complex structures $J$ is rich enough to allow deformations of  pseudoholomorphic curves.  These results should be known to experts in the field, see \cite{Bi} or \cite{U}.

One purpose of this note is to introduce the relative symplectic
cone and prove a version of the gluing formula for fiber sums along
$T^2$.  Missing from the examples in Section \ref{relcone} is the K3 surface.  This will be detailed in a further paper, thus rounding out the known examples of symplectic manifolds with Kodaira dimension 0.

The second author is supported by the NSF grant 0604748.

\section{\label{relcone}The relative symplectic cone}

Let $M$ be an oriented  manifold and $V$ an oriented
codimension 2 submanifold, not necessarily connected.  Throughout this section, it will be necessary to carefully distinguish the class of $V$, denoted $[V]\in H_2(M)$, and the specific submanifold $V$.  Denote the
Poincar\'e dual of $[V]$ by $[V]^D$.  We make the following definition:

\begin{definition}A relative symplectic form on the pair $(M, V)$ is
 an orientation compatible symplectic form on $M$ such that
$\omega\vert_V$ is an  orientation compatible symplectic form on $V$.
The relative symplectic cone of $(M, V)$ is
 \begin{equation}
{\mathcal C}_M^V=\{\alpha\in H^2(M)|\;
[\omega]=\alpha,\;\omega \mbox{ is a relative symplectic form on $(M, V)$}\}.
\end{equation}
\end{definition}

The following lemma follows directly from the definition of the relative symplectic cone:

\begin{lemma}
\label{incl}${\mathcal C}_M^V\subset \{\alpha\in{\mathcal C}_M\;\vert\;\alpha\cdot [V]^D>0\}\subset \{\alpha\in{\mathcal P}_M\;\vert\;\alpha\cdot [V]^D>0\}$
\end{lemma}

Obviously there are maps
\begin{equation}
\label{res}{\mathcal C}_M^V\hookrightarrow {\mathcal C}_M, \quad {\mathcal C}_M^V\rightarrow {\mathcal C}_V.
\end{equation}
In fact, if $V$ is the disjoint union of $V_0$ and $V_1$, then there
is also a  map ${\mathcal C}_M^V\rightarrow {\mathcal C}_M^{V_i}$.  Note that the restriction mapping ${\mathcal C}_M^V\rightarrow {\mathcal C}_V$ is by no means generically injective.  The following fact relates the relative cones to the symplectic cone:

\begin{lemma}\label{union}Let $\mathcal V$ denote the set of oriented codimension 2 submanifolds of $M$.  Then
\[\bigcup_{V\in\mathcal V}\mathcal C_M^V=\mathcal C_M\]
\end{lemma}

\begin{proof} The inclusion $\bigcup_{V\in\mathcal V}\mathcal C_M^V\subset\mathcal C_M$ follows from \ref{res}.  Consider now a symplectic class $\omega\in \mathcal C_M$.  We distinguish two cases: $b^+=1$ and $b^+>1$.  If $b^+=1$, then Prop 4.2 or Prop 4.3 (together with Lemma 3.5) in \cite{TJLL} shows that the class $k\omega$ for $k$ large enough is represented by a $\omega$-symplectic surface.  Hence, if $V$ represents the class $k\omega$, then $\omega\in \mathcal C_M^V$.

If $b^+>1$, then the canonical class of $(M,\omega)$ for some almost complex structure $J$ taming $\omega$ is represented by a $\omega$-symplectic surface, see \cite{T}.  Hence, if $V$ represents the canonical class $-c_1(M,\omega)=-c_1(M,J)$, then $\omega\in\mathcal C_M^V$. 
\end{proof}

The proof shows, that if $b^+>1$, we need only consider submanifolds $V$ which are representatives of a canonical class $K_\omega$ of $M$ if we wish to understand $\mathcal C_M$ with respect to the relative cone.  In particular, this shows that $\mathcal C_M\subset {\mathcal P}^{c_1(M)}\cup{\mathcal P}^{-c_1(M)}$ if $b^+>1$ and $M$ is minimal K\"ahler, which is of interest in connection with Conjecture \ref{conj}.  Furthermore, it seems natural to wonder, whether there exist a finite set of submanifolds $\mathcal V_f\subset \mathcal V$, such that they completely determine the symplectic cone of $M$.  With respect to this question, a key observation connecting the relative symplectic cone and the symplectic cone is

\begin{lemma}\label{VD}  Assume $c_1(M,\omega)=a[V]^D$ with $a\ne 0$ for some symplectic form $\omega$ and
\begin{equation}
{\mathcal C}^V_M=\{\alpha\in {\mathcal
P}\;\vert\;\alpha\cdot [V]^D> 0\},
\end{equation}
then
 ${\mathcal
P}^{c_1(M,\omega)}\cup{\mathcal P}^{-c_1(M,\omega)}\subset{\mathcal C}_M$.  If $M$ is minimal K\"ahler, then ${\mathcal C}_M={\mathcal
P}^{c_1(M)}\cup{\mathcal P}^{-c_1(M)}$
\end{lemma}

\begin{proof}Clearly, if $\alpha\in {\mathcal C}^V_M$, then $\alpha\cdot c_1(M,\omega)\ne 0$.
Furthermore, the definition of relative symplectic form and relative symplectic cone assumes an orientation
of the submanifold $V$. We could however use the opposite orientation on
$V$ as well, thus changing the sign $\alpha\cdot [V]^D<0$ and hence the sign of $\alpha\cdot c_1(M,\omega)$. Therefore, ${\mathcal
P}^{c_1(M,\omega)}\cup{\mathcal P}^{-c_1(M,\omega)}\subset {\mathcal C}_M$ by Lemma \ref{union}.  If $M$ is minimal K\"ahler, the result follows from the inclusion ${\mathcal C}_M\subset{\mathcal
P}^{c_1(M)}\cup{\mathcal P}^{-c_1(M)}$ (see \cite{TJL1}). 
\end{proof}

We now proceed to calculate the relative cone for certain submanifolds $V$ for two classes of symplectic manifolds:  $T^2$ bundles over $T^2$  and manifolds with $b^+=1$.

\subsection{$T^2$-bundles over $T^2$}\label{Torus}  The total spaces $M$ of such bundles have been studied and classified by Sakamoto-Fukuhara \cite{SF}, Ue \cite{Ue} and Geiges \cite{Ge}.  In particular, with one exception, they all admit symplectic structures compatible with the bundle structure; in the case of a primary Kodaira surface this bundle structure must be specified as it is not unique.  Moreover, the relative symplectic cone with respect to the fiber torus $T^2_f$ has been determined explicitly by Geiges.  

 In \cite{SF}, the manifolds $M$ are classified according to the monodromy $A,B$ of the bundle and the Euler class $(x,y)$.  A manifold $M$ is determined by the tuple $(A,B,(x,y))$.  In \cite{Ue}, the total spaces are classified according to their geometric type as defined by Thurston.  Furthermore, an explicit representation of each is given in terms of generators of $\Gamma$ such that $M={\mathbb R}^4\backslash \Gamma$.  For example, the four torus $T^4$ is given by the following data:  $(Id,Id,(0,0))$ (Id is the 2x2 identity matrix) with geometric type $E^4$ and $\Gamma=\mathbb Z^4$, i.e. $T^4=\mathbb R^4\backslash \mathbb Z^4$.  From the explicit presentation of the generators of $\Gamma$, Geiges constructs symplectic forms, thereby determining the symplectic cones as well as the relative cones with respect to the fiber torus $T^2_f$.  In the following we denote the class of the fiber torus $T_f^2$ by $F$.  We collect the data in the following table, details can be found in \cite{Ue} and \cite{Ge}: 
\[
\begin{array}{|c|c|c|c|}\hline
\mbox{type}&b_1&\mathcal C_M&\mathcal C_M^{T_f^2}\\\hline
T^4&4&\mathcal P_M&\mathcal P_M^F\\\hline
 \mbox{primary Kodaira surface}&3&\mathcal P_M & \mathcal P_M^F\\\hline
\mbox{hyperelliptic surface}&2&\mathcal P_M&\mathcal P_M\\\hline
(d)&2&\mathcal P_M&\emptyset\\\hline
(e)-(h)&2&\mathcal P_M&\mathcal P_M\\\hline
\end{array}
\]

Note that the class of $T^2$-fibrations over $T^2$ provides a full range of possible relative cones, from $\emptyset$ to the maximal possible cone, see Lemma \ref{incl}.

\subsection{Manifolds with $b^+=1$}
The symplectic cone in this case is determined in \cite{TJLL}:

\begin{theorem}\label{b1}Let $M$ be a closed, oriented 4-manifold with $b^+=1$ and $\mathcal C_M$ nonempty.  Let $\mathcal E$ denote the set of all exceptional classes of $M$.  Then
\[
\mathcal C_M=\{e\in\mathcal P_M\;\vert\; 0<\vert e\cdot E\vert\;\;\mbox{for all }\;E\in \mathcal E\}.
\]
\end{theorem}

In this section, we will determine the relative cone with respect to a submanifold $V$ for manifolds with $b^+=1$.  This result depends in large part upon the results in the Appendix as well as on results found in \cite{TJLL}.  A trivial corollary following from Thm. \ref{b1} and Lemma \ref{incl} is

\begin{cor}\label{trivial inclusion}
\[
\mathcal C^V_M\subset \{e\in\mathcal P_M\;\vert\; 0<\vert e\cdot E\vert\;\;\mbox{for all }\;E\in \mathcal E\mbox{ and }e\cdot [V]^D>0\}
\]

\end{cor}

The following result shows that this inclusion is actually an equality: 

\begin{theorem} \label{li1}Let $M$ be 4--manifold with  $b^+=1$ and
nonempty $C_M$. If $V$ is an oriented  submanifold for which
${\mathcal C}_M^V\ne \emptyset$,  then
\[
\mathcal C^V_M= \{e\in\mathcal P_M\;\vert\; 0<\vert e\cdot E\vert\;\;\mbox{for all }\;E\in \mathcal E\mbox{ and }e\cdot [V]^D>0\}=\mathcal P_{\mathcal E}^{[V]}.
\]
In particular, if $M$ is minimal, then
$${\mathcal C}_M^V=\mathcal P^{[V]}=\{e\in \mathcal P_M|e([V])>0\}.$$
\end{theorem}

\begin{proof}  The following Lemma will be central to the proof.  We will complete the proof of the Theorem before proving the Lemma.  The proof of the Theorem relies upon the relative inflation procedure in \cite{Bi}.  For this method to work, we need to find two symplectic submanifolds satisfying certain restrictions.  One of these will be $V$, the other will be constructed in the Lemma.  We begin with a definition:

\begin{definition}Fix $\omega\in \mathcal C_M$.  Let $\mathcal E_\omega$ be the set of exceptional curves $E\in\mathcal E$ which can be represented by an $\omega$-symplectic sphere of self-intersection -1.
\end{definition}

 \begin{lemma} \label{rel submanifold} Let us fix a relative symplectic form $\omega$ on $(M, V)$.
 For any $A\in H_2(M;\mathbb Z)$ with
\[A\cdot E >0\mbox{ for all } E\in \mathcal E_\omega,\]
 $$A\cdot
 [V]>0,\hspace{3mm} A\cdot A>0,\hspace{3mm} A[\omega]>0$$
  $$(A-K_{\omega})\cdot [\omega]>0,\hspace{3mm} (A-K_{\omega})\cdot
  (A-K_{\omega})>0,\hspace{3mm} (A-K_{\omega})\cdot [V]>0,$$
there exists an $\omega-$symplectic submanifold in the class $A$,
intersecting $V$ transversally and positively.
\end{lemma}

It follows from Corollary \ref{trivial inclusion}, that we need only to show ${\mathcal P}_{\mathcal E}^{[V]}\subset {\mathcal C}_M^V$.  We will do this in two steps:  First we prepare every class in  ${\mathcal P}_{\mathcal E}^{[V]}$ to satisfy the assumptions in Lemma \ref{rel submanifold}.  In the second step, we will use the relative inflation process in \cite{Bi} to
construct relative symplectic forms for the pair $(M, V)$ in the class
$[\omega]+tA$ for any $t>0$ for suitably chosen $A$ and any $\omega\in {\mathcal C}_M^V$.  Fix a class $\omega\in \mathcal C_M^V$, we can assume that $[\omega]$ is an integral class.

Let $e\in {\mathcal P}_{\mathcal E}^{[V]}$ be any class such that $e\cdot [\omega]>0$.

{\bf 1.1 $M$ is not rational or ruled.} In this case, \cite{TJL2} shows that $\mathcal E=\{\pm F_1,...,\pm F_l\}$ for some finite $l$.  Furthermore, \cite{TJLL} shows that any symplectic canonical class of $M$ is of the form $K=\pm V\pm F_1\pm...\pm F_l$ where $V$ is, up to sign, the unique symplectic canonical class of the minimal model.

Let $e\in {\mathcal P}_{\mathcal E}^{[V]}$, then results in \cite{TJL2} and \cite{TJLL} show that there exists a choice of canonical class $K$ associated to $\omega$ such that the following properties hold:
\begin{enumerate}
\item $e\cdot [V]>0$ and
\item $e\cdot E>0$ for all $E\in \mathcal E_\omega$.
\end{enumerate}

{\bf 1.2 $M$ is rational or ruled.} Results in \cite{TJBL}, show that it is possible to reduce every class $e$ with positive square to a special reduce class via a diffeomorphism $\phi$.  This diffeomorphism takes $\mathcal C_M^V$ to $\phi^*(\mathcal C_M^V)=\mathcal C_M^{\phi V}$.  In particular, it preserves the inequality $e\cdot [V]>0$.  Moreover, it is shown in \cite{TJBL} and \cite{TJLL}, that $\phi^*e\cdot E>0$ for all $E\in \mathcal E_\omega$.

{\bf 2. Relative Inflation.}  Let $e$ be one of the classes obtained above, i.e. $e$ in the first and $\phi^*e$ in the second.  Apply Lemma \ref{rel
submanifold} to the class $A=le-[\omega]$
for $l>>0$ and use the relative inflation procedure in \cite{Bi} with $Z=V$ and $C=$ symplectic surface in class $A$ obtained by Lemma \ref{rel submanifold}. This proves that $le$ is in the relative cone. Hence $e$ is
itself in the relative cone. It also follows that any real multiple
of an integral class $e$ in $\mathcal P_{\mathcal E}^{[V]}$ is in the relative
cone ${\mathcal C}_M^V$.

Observe that $\mathcal P_{\mathcal E}^{[V]}$ is an open convex cone.  Therefore,
for any $\alpha$ in $\mathcal P_{\mathcal E}^{[V]}$, we can write
$\alpha=\sum_{i=1}^p\alpha_i$, where the rays of $\alpha_i$ are in
$\mathcal P_{\mathcal E}^{[V]}$, arbitrarily close to that of $\alpha$, and each
$\alpha_i=s_i\beta_i$ for some positive real number $s_i$ and an
integral class $\beta_i$.  Note that $\beta_i\cdot \beta_j >0$ for all $i,j$.  Inductively it can be shown that for each
$q\leq p$, $\sum_{i=1}^q\alpha_i$ is in the relative cone ${\mathcal
C}_M^V$:

First we choose a relative symplectic form $\omega_1$ for the pair $(M, V)$
with $[\omega_1]=\alpha_1$. For a large integer $l$, since
$$[\beta_2]\cdot [\beta_2]>0, \beta_2\cdot [\omega_1]>0, \beta_2\cdot [V]>0,$$ we can apply Lemma
\ref{rel submanifold} to $A=l[\beta_2]$. By relative inflation, we
find that $\alpha_1+\alpha_2=\alpha_1+\frac{s_2}{l}A$ is in the
relative cone.

Now choose a symplectic form $\omega_2$ for the pair $(M, V)$ with
$[\omega_2]=\alpha_1+\alpha_2$.  This completes the argument. 

\end{proof}

\begin{proof}(of Lemma \ref{rel submanifold})  This relies on Proposition 4.3 in \cite{TJLL} and the genericity results of the Appendix.

Let us first recall the notion of $A$ being $J-$effective and simple
considered in \cite{Bi}.  $A$ is $J-$effective and simple if, for a
generic choice of $k(A)=\frac{-K_{\omega}(A)+A\cdot A}{2}\geq 0$
distinct points $\Omega_{k(A)}$ in $M$, there exists a connected $J$-holomorphic
submanifold $C\subset M$ which represents $A$ and passes through all
the $k(A)$ points.

Proposition 4.3 in \cite{TJLL} implies that any $A$ with
\[A\cdot E\ge 0\mbox{ for all }E\in \mathcal E,\]
$$A\cdot A>0,\hspace{3mm} A\cdot [\omega]>0,$$
$$(A-K_{\omega})\cdot (A-K_{\omega})>0,\hspace{3mm} (A-K_{\omega})\cdot
[\omega]>0,$$ is $J-$effective and simple for a generic $J$ tamed by
$\omega$.  More precisely,  $A$ is $J$-effective for all such $J$, i.e. there exists a connected $J$-holomorphic curve $\phi:\Sigma\rightarrow M$ representing $A$ such that $Im(\phi)$ is a submanifold passing through all $k(A)$ points.   However $\phi$ may be a multiple cover or have multiple components.

Let $\mathcal J_V$ be the space of
$\omega-$tamed almost complex structures making $V$ a
pseudo-holomorphic submanifold.  Then for every $J\in\mathcal J_V$ there exists a $J$-holomorphic curve $C$ in class $A$ by the previous considerations.  We consider connected nodal $J$-holomorphic curves $C$ representing $A$ with multiple components each having as their image a $\omega$-symplectic submanifold.  However, for the purposes of this Lemma,  we need to exclude components which lie in $V$.  Therefore consider a connected $J$-holomorphic curve $C$ representing $A$:  The curve $C$ is given by a collection $\{(\phi_i,\Sigma_i)\}$ of maps $\phi_i$  and Riemann surfaces $\Sigma_i$.  We want to reformulate this as a collection $ \{(\varphi_i,C_i,m_i)\}$ of simple maps, submanifolds and multiplicities.  If $\phi_i$ is a multiple cover, we replace it by a simple map $\varphi_i$ with the same image and an integer $m_i$ tracking the multiplicity.  We also combine maps which have same image, adding the multiplicities and keeping only one copy of the map.  Ultimately, we replace all the maps $\phi_i$ by simple embeddings $\varphi_i$ with image $C_i$.  We can therefore decompose the class $A$ as $\sum_im_i[C_i]$.  Note that we are only interested in the submanifolds $C_i$, so we do allow the class $[C_i]$ to be divisible.  However, we want the class to represent a submanifold, hence correspond to an embedding $\varphi_i$.  We denote
\begin{enumerate}
\item components with $[C_i]^2<0$ by $B_i$ and
\item those with $[C_i]^2\ge 0$ which do not lie in $V$ by $A_i$.
\end{enumerate}
The class $A_i$ could be a multiple of the class $[V]$, however due to our decomposition above we consider only maps which are not multiple covers of $V$.  Furthermore, we could have a component $[C_0]=[V]$, with a multiplicity $m\ge 1$, which is a (multiple) cover of $V$.  Thus, $A=m[V]+\sum_i m_iB_i+r_iA_i$.

We begin with the negative square components:  Let $B\cdot B<0$.  Consider the moduli space of $\mathcal M(B, J, g)$ of pairs $(u,
j)$, where $j\in \mathcal T_g$, the Teichm\"uller space of a closed
oriented surface $\Sigma_g$ of genus $g$, and $u:(\Sigma_g, j)\to
(M, J)$ is a somewhere injective $(j, J)-$holomorphic curve in the
class $B$. If $B\ne [V]$, then arguments similar to those in the proof of Lemma \ref{genericA} show that
for generic $J\in \mathcal J_V$, the spaces $\mathcal M(B, J, g)$
are smooth manifolds of dimension $2(-K_{\omega}\cdot B+g-1)+\dim
G_g$. 

By the adjunction formula the space of non-parametrized $J-$holomorphic curves
has dimension
$$2(-K_{\omega}\cdot B+g-1)\leq 2B\cdot B-(2g-2).$$
Thus if $B\cdot B<0$, $\mathcal M(B, J, g)$ is non-empty only if
$g=0$ and $B\cdot B=-1$.  We conclude,  that for  a generic $J\in \mathcal J_V$, the only irreducible
components of a cusp $A-$curve with  negative self-intersection (except possibly the component $C_0$) have $B^2=-1$.  In particular, note that all $B_i\in \mathcal E_\omega$ and $k(B_i)=0$.  Hence, our assumptions imply that $A\cdot B_i > 0$ for all $i$.   

Now let us divide the proof into two cases, in both we shall use the genericity results proven in the Appendix:

{\bf Case 1. $k([V])\ge 0$:}  The condition $k([V])\ge 0$ can be rewritten using the adjunction formula to state that $[V]^2\ge0 $ if $g([V])\ge 1$ and $[V]^2\ge -1$ at worst, if $g([V])=0$.  In the following we will allow the case $A=m[V]$.

The results of the Appendix, in particular Lemma \ref{genericA} and Lemma \ref{jV}, show, that we can find a generic set of pairs $(J,\Omega_{k(A)})$ such that $k(A_i)\ge 0$, each curve in class $A_i$ resp. $[V]$ meets at most $k(A_i)$ resp. $k([V])$ generic points and $\sum k(A_i)+k([V])\ge k(A)$.  

Further, if $k(A_i)\ge 0$ and $A_i^2\ge 0$, then $k(r_iA_i)\ge 0$ for any positive integer $r_i$:
\[
0\le 2k(A_i)\le 2r_ik(A_i)=-K_{\omega}(r_iA_i)+r_iA_i\cdot A_i\le
\]
\[
\le -K_{\omega}(r_iA_i)+r_i^2A_i\cdot A_i=2k(r_iA_i).
\]
Note that this holds in particular for $m[V]$.

For such a generic choice of $(J,\Omega_{k(A)})$, let $C$ be a connected curve representing $A$, which contains the $k(A)$ distinct points of $\Omega_{k(A)}$.  Then
\[
2k(A)=-K_{\omega}(m[V])+\sum_i-K_\omega(m_iB_i)+\sum_i-K_{\omega}(r_iA_i)+
\]
\[
+m^2[V]^2+\sum_im_i^2B_i^2+\sum_ir_i^2A_i^2+2\sum_im[V]m_iB_i+2\sum_im[V]r_iA_i+
\]
\[
+2\sum_{i> j}m_im_jB_i\cdot B_j+2\sum_{i> j}m_ir_jB_i\cdot A_j+2\sum_{i> j}r_ir_jA_i\cdot A_j\ge
\]
\[
\ge 2mk([V])+2\sum_i r_ik(A_i)+ (m^2-m)[V]^2+2\sum_im[V]r_iA_i+2\sum_{i> j}r_ir_jA_i\cdot A_j+
\]
\[
+\sum_i (m_i^2-m_i)B_i^2+2\sum_{i> j}m_im_jB_i\cdot B_j+2\sum_{i> j}m_ir_jB_i\cdot A_j+2\sum_im[V]m_iB_i
\]
If $[V]^2=-1$, then denote $B_0=[V]$ and include it in the following estimate.  Fix an $i$ and consider the terms in the last line.  They can be rewritten as 
\[
2m_iA\cdot B_i-2m_i^2B_i^2+(m_i^2-m_i)B_i^2\;=\;2m_iA\cdot B_i+m_i^2+m_i\;\ge 0
\]
and thus we obtain the estimate
\[
2k(A)\ge 2k([V])+2\sum_i k(A_i).
\]
Hence either $k(A)> k([V])+\sum_i k(A_i)$ or  the following hold:
\begin{itemize}
\item $m_i=0$ for all $i$, i.e. there are no components of negative square,
\item $A_i\cdot A_j=0=A_i\cdot [V]$ for $i\ne j$,
\item if $[V]^2\ge 0$, then $m=1$ or $k([V])=0$ and $[V]^2=0$ and
\item $r_i=1$ or $k(A_i)=0$.
\end{itemize}
Therefore, the curve $C$ representing $A$ is an embedded $J$-holomorphic submanifold with a single non-multiply covered component containing the set $\Omega_{k(A)}$ with $J\in \mathcal J_V$ or $k(A)=0$.  

In the following cases we are done:
\begin{enumerate}
\item $A\ne m[V]$
\item $A=m[V]$ and $k([V])>0$: The results above imply that $m=1$.  Choose $\Omega_{k(A)}$ such that it contains a point not in $V$.  Then $C$ does not lie in $V$ and intersects $V$ locally positively.
\end{enumerate}
We are left with the following case: $A=m[V]$ and $k([V])=0$.  However, in this case the previous results show that either $m=1$ or $[V]^2=0$.  The latter is excluded by the assumption $A^2>0$.  The former would mean $A=[V]$ and thus
\[
0=2k([V])=[V]^2-K_\omega [V]=[V](A-K_\omega) 
\]   
which is also excluded by assumption.

In all cases, we can perturb $C$ to be transverse to $V$, see \cite{M4}, \cite{M5}.

{\bf Case 2: $k([V])<0$:} In this case, the results of the Appendix show, that  we can find a generic set of almost complex structures, such that $V$ is rigid and there are no other curves in class $[V]$.  In the following, we choose only complex structures from this set.

Even though we are working in the case $k([V])<0$, it is possible for a multiple class $m[V]$ to have $k(m[V])\ge 0$.  For this reason, we will distinguish the following two objects:
\begin{enumerate}
\item Classes $A_i=m_i[V]$ which correspond to components of the curve $C$ in class $m[V]$, but which are NOT multiple covers of a submanifold in class $[V]$.  If $[V]^2< 0$, then positivity of intersections shows that any curve $C$ can contain at most one component in class $m[V]$  for all $m$ and this component must coincide with the manifold $V$.  This situation was studied in greater generality in \cite{Bi}.  Furthermore, if  $[V]^2\ge 0$ and a class $A_i=m[V]$ occurs, then the results of Lemma \ref{genericA} apply.  We may therefore assume , that $A_i^2\ge 0$ in the following.
\item The specific "class" $mV$ which corresponds to components which have as their image the submanifold $V$.
\end{enumerate}
Note further, that we can choose our almost complex structures such that the components corresponding to $mV$ are rigid, while those in $m_i[V]$ are not.  Such a decomposition is not necessary in the case $k([V])\ge 0$, as we can choose a generic set of pairs such that $\mathcal K_V^{[V]}(J,\Omega)$ is smooth (see the Appendix for details), however we do not know if $V$ is an element in this set, nor does this matter for the calculation.  In the current situation, the specific submanifold $V$ acts differently than other elements in the class $[V]$.

We now proceed as in \cite{Bi}: Consider the class $\tilde A=A-mV=\sum r_iA_i$.  We assume that such a decomposition is possible, i.e. there exists a not necessarily connected pseudoholomorphic curve in class $\tilde A$.  The case $A=mV$ will be considered afterwards.  We need to show, that there exists a generic set of pairs $(J,\Omega_{k(A)})$ such that $k(A)>\sum k(A_i)$.  Proceeding exactly as in \cite{Bi}, we obtain the following two estimates: $\sum k(A_i)\le k'(\tilde A)$, where $k'$ is the modified count defined by McDuff \cite{M3}.  Furthermore
\[
2k(A)-2k'(\tilde A)= m \underbrace{(A-K_\omega)\cdot V}_{>0 \mbox{ by assumption}}+ \mbox{non-negative terms}
\]
and hence, combining all these inequalities,  for generic pairs $(J,\Omega_{k(A)})$ we obtain $k(A)>\sum k(A_i)$.  We therefore conclude, that we can rule out such a decomposition, hence $"A=mV"$ or $A=\sum r_iA_i$.  In the latter case we are done by the same line of argument as in the $k([V])\ge 0$ case, albeit with the added restriction on the almost complex structures that $V$ is rigid and no further curves in class $[V]$ exist.  The case $"A=mV"$ corresponds to the case $A=m[V]$ and $k(A)\ge 0$, more precisely,
\[
0<m(A-K_\omega)\cdot V=mA\cdot V-K_\omega\cdot V=A^2-K_\omega\cdot A=k(A).
\]
Thus, applying  Lemma \ref{genericA}, we can find a generic set of pairs $(J,\Omega_{k(A)})$ such that $A=m[V]$ is represented by an embedded curve with deformations.  Choosing $\Omega_{k(A)}$ such that it contains a point not in $V$ ensures that a representative of $A$ in this case does not lie in $V$.

As before, we can make the curve $C$ transverse to $V$.

 \end{proof}

\section{The gluing formula \label{cutsum}}

We now return to the question posed in the introduction. In this
section we provide the theoretical framework necessary to answer
this question for so called "good" sums. We first review the
symplectic sum and cut operations. This leads to the definition of a
good sum and the subsequent homological reformulation of these
operations.

\subsection{Smooth fiber sum}

Let $X$, $Y$ be $2n$-dimensional smooth
manifolds. Suppose we are given
codimension 2 embeddings $j_*:V\rightarrow *$ into $X$ and $Y$ of  a
smooth closed oriented manifold $V$ with normal bundles $N_*V$.
Assume that the Euler classes of the normal bundle of the embedding
of $V$ in $X$ resp. $Y$ satisfy $e(N_XV)+e(N_YV)=0$ and fix a
fiber-orientation reversing bundle isomorphism $\Theta:
N_XV\rightarrow N_YV$.  By canonically identifying the normal
bundles with a tubular neighborhood $\nu_*$ of $j_*(V)$, we obtain
an orientation preserving diffeomorphism $\varphi: \nu_X\backslash
j_X(V)\rightarrow \nu_Y\backslash j_Y(V)$ by composing $\Theta$ with
the diffeomorphism that turns each punctured fiber inside out.  This
defines a gluing of $X$ to $Y$ along the embeddings of $V$ denoted
$M=X\#_{(V,\varphi)}Y$.  The diffeomorphism type of this manifold is
determined by the embeddings $(j_X,j_Y)$ and the map $\Theta$.  Note also, that if $V$ has trivial normal bundle, then this construction should actually be viewed as a sum along $V\times S^1$.

In the rest of the paper, whenever we consider a fiber sum, we fix
$V$, the embeddings $j_*$ and the bundle isomorphism $\Theta$
without necessarily explicitly denoting either.  This fixes the
homology of the manifold $M=X\#_{(V,\varphi)}Y$.

\begin{example} \label{T4T4} Consider for example the torus $T^4=T_f^2\times T^2$, where the first
factor is the fiber direction.  Let  $M=T^4\#_{T_f^2}T^4$ be the sum
along the fiber $T_f^2$. Then $M$ is actually $T_f^2\times \Sigma_2$, as can be seen from the following:

\[\begin{diagram}
\node{T_f^2}\arrow{e}\node{T^4}\arrow{s}\node{T^4}\arrow{s}\node{T_f^2}\arrow{w}\\
\node[2]{T^2}\node{T^2}
\end{diagram}
\Longrightarrow
\begin{diagram}
\node{T_f^2}\arrow{e} \node{T^4\#_{T_f^2}T^4=M}\arrow{s}\\
\node[2]{T^2\#T^2=\Sigma_2}
\end{diagram}\]

\end{example}

\subsection{Symplectic sum and symplectic cut}

We briefly describe the symplectic sum construction $M=X\#_VY$  as
defined by Gompf \cite{Go2} (see also McCarthy-Wolfson \cite{MW}).  Assume $X$ and $Y$ admit symplectic forms $\omega_X,\omega_Y$ resp.  If the embeddings $j_*$ are symplectic with respect to these forms, then
we obtain $M=X\#_{(V,\varphi)}Y$ together with a symplectic form
$\omega$ created from $\omega_X$ and $\omega_Y$.  Furthermore, it
was shown in \cite{Go2} that this can be done without loss of
symplectic volume.

Furthermore, Gompf showed that the symplectic form $\omega$ thus
constructed on $M=X\#_{(V,\varphi)}Y$ from $\omega_X,\omega_Y$ is
unique up to isotopy.  This result allows one to construct a smooth
family of isotopic symplectic sums $M=X\#_{(V,\varphi_\lambda)}Y$
parametrized by $\lambda \in D^2\backslash\{0\}$ as deformations
with a singular fiber $X\sqcup_VY$ over $\lambda=0$ (see
\cite{IP5}).  Therefore, we suppress $\varphi_\lambda$ from the
notation, choosing instead to work with an isotopy class where
necessary.

Thus, a symplectic sum will be denoted by $M=X\#_VY$,  a symplectic
class $\omega$ on the sum will denote an isotopy class.

The symplectic cut operation of Lerman \cite{Le} functions as
follows: Consider a symplectic manifold $(M,\omega)$ with a
Hamiltonian circle action and a corresponding moment map
$\mu:M\rightarrow {\bf \mathbb R}$.  We can assume that 0 is a regular value, if
necessary by adding a constant.  We can thus cut $M$ along
$\mu^{-1}(0)$  into two manifolds $M_{\mu>0}$ and $M_{\mu<0}$, both
of which have boundary $\mu^{-1}(0)$. If we collapse the
$S^1$-action on the boundary, we obtain manifolds
$\overline{M_{\mu>0}}$ and $\overline{M_{\mu<0}}$ which contain  a
real codimension 2 submanifold $V=\mu^{-1}(0)/S^1$.  If we
symplectically glue $\overline{M_{\mu>0}}$ and
$\overline{M_{\mu<0}}$ along $V$ we obtain again $M$.

Note that the above construction is local in nature,
thus the moment map and the $S^1$ action need only be defined in a neighborhood of the cut.

The symplectic structure $\omega$ restricted to $M_{\mu>0}$ and
$M_{\mu<0}$ reduces to a symplectic structure on
$\overline{M_{\mu>0}}$ and $\overline{M_{\mu<0}}$ which have the
same value on $V$.  This motivates the sum decomposition of the symplectic cones in section \ref{main}.

A symplectic cut is only possible on a symplectic manifold, thus, when  discussing a symplectic cut on $M=X\#_VY$, we implicitly consider only those isotopy classes allowing moment maps $\mu$ with $V=\mu^{-1}(0)\backslash S^1$.

In the case $M=T^4\#_{T^2_f}T^4$ it is possible to understand the geometric construction underlying the cut:  Consider $M=T^4\#_{T^2_f}T^4=T^2\times \Sigma_2$ and view $\Sigma_2$ such that we have the holes on either end and a cylindrical connecting piece $S$ in between.  Furthermore, in $M$ this copy of $\Sigma_2$ is transverse to the fiber $T^2$.  Choose local coordinates $(\lambda,\theta,t)$ on $S\times T^2$, $(\lambda,\theta)\in [-1,1]\times [0,2\pi]$ coordinates on $S$, $t$ a coordinate on $T^2$.  Consider an $S^1$ action on the second coordinate stemming from the Hamiltonian $\mu:S\rightarrow {\bf \mathbb R}$ given by $\mu(\lambda,\theta,t)=\lambda$.  Locally, any symplectic form is given by $\omega=a \;d\lambda\wedge d\theta + b\;dt + \Omega$, with $\Omega(\frac{\partial}{\partial \theta},\cdot)=0$ and $a\in {\bf \mathbb R}$ nonzero.

The symplectic cut defined by $\mu$ produces $M_{\mu<0}= T^2_b\times T^2$ with $T^2_b$ a punctured torus with boundary $S^1$.  Over each point of the boundary, there is a fiber $T_f^2$, hence $\partial M_{\mu<0}=T^3$.  Collapsing this boundary under the $S^1$ action produces $T^4=T_f^2\times T^2$.  In particular, the action maps $d\lambda\wedge d\theta$ to local coordinates on a neighborhood of the collapsed boundary on $T^2_b$ without loss of volume.

\subsection{\label{main}${\mathcal C}^V_{X\#_VY}$ from ${\mathcal C}_X^V$ and ${\mathcal C}_Y^V$}

\subsubsection{The cone of sum forms}
 We are interested in the symplectic cone ${\mathcal C}_M$ of the
4-manifold $M$.  Suppose this manifold can be obtained as a
symplectic sum $M=X\#_VY$. Let us fix the symplectic embeddings as
well as the map $\Theta$.  In the following, we will distinguish
between the manifold $M$ and the specific viewpoint as a symplectic
sum from $X$ and $Y$ along $V$ by explicitly denoting $M=X\#_VY$.
Accordingly, we define the following symplectic cone associated to
the symplectic sum:

\begin{definition}Suppose that $M=X\#_VY$. Define the cone of sum forms, ${\mathcal
C}_{X\#_VY}^{sum}$,  to be the set of  classes of symplectic
forms on $M$ which can be obtained by summing $X$ and $Y$ with
symplectic embeddings $\tilde j_*$ and bundle map $\tilde\Theta$
isotopic to the fixed choice $j_*$ and $\Theta$.
\end{definition}

We obtain the following result:

\begin{theorem}\label{t1} For a symplectic manifold $M=X\#_VY$,
\begin{equation}
{\mathcal C}_{X\#_VY}^{sum}=\Psi(\Phi^{-1}({\mathcal C}_X^V\oplus {\mathcal C}_Y^V))
\end{equation}
where  $\Phi,\Psi$ are the maps on cohomology corresponding to the inclusion of $X,Y$ into $X\sqcup_VY$ and the projection of $X\#_VY$ onto the singular manifold $X\sqcup_VY$ respectively (see (\ref{pd}) below).
\end{theorem}

\begin{proof}Consider the following maps on cohomology:
\begin{equation}
\begin{diagram}\label{pd}
\node{H^2(X\sqcup_VY)}\arrow{e,t}{\Psi}\arrow{s,l}{\Phi}\node{H^2(X\#_VY)}\\
\node{H^2(X)\oplus H^2(Y)}\\
\end{diagram}
\end{equation}
Let  $C_{X\sqcup_VY}:=\Phi^{-1}({\mathcal C}_X^V\oplus {\mathcal C}_Y^V)\in
H^2({X\sqcup_VY})$. We can view this set as the collection of classes of
symplectic forms in $X$ and $Y$ which are symplectic and equal on
$V$.  More precisely, we obtain $C_{X\sqcup_V Y}$ by pulling back of ${\mathcal C}_V$ under the restriction maps $r_X,r_Y$ from (\ref{res}):
\[
\begin{diagram}
\node{C_{X\sqcup_VY}}\arrow{e}\arrow{s}\node{{\mathcal C}_X^V}\arrow{s,l}{r_X}\\
\node{{\mathcal C}_Y^V}\arrow{e,t}{r_Y}\node{{\mathcal C}_V}\\
\end{diagram}
\]
The symplectic sum takes $(X,\omega_X)$ and $(Y,\omega_Y)$ and
produces a symplectic manifold $(X\#_VY,\omega)$.  This will work
for any relative symplectic forms $\omega_X$ and $\omega_Y$
identical on $V$ (By identical we mean that the symplectomorphism
used to produce the symplectic singular manifold $X\sqcup_VY$ maps
these two forms symplectically to each other along $V$.).  Thus any
symplectic class $(\alpha_X,\alpha_Y)\in C_{X\sqcup_VY}$ can be
summed to produce a symplectic class $\alpha \in {\mathcal
C}_{X\#_VY}^{sum}$.  Therefore $\Psi(\Phi^{-1}({\mathcal C}_X^V\oplus
{\mathcal C}_Y^V))\subset {\mathcal C}^{sum}_{X\#_VY}$.

On the other hand, given any  symplectic class in ${\mathcal
C}_{X\#_VY}^{sum}$,  any symplectic representative $\omega$ of such
a class can be symplectically cut, such that the manifolds $(X,V)$
and $(Y,V)$ result with symplectic forms $\omega_X$ and $\omega_Y$
which agree on $V$.  Hence, $\Psi^{-1} {\mathcal
C}_{X\#_VY}^{sum}\subset\Phi^{-1}({\mathcal C}_X^V\oplus {\mathcal C}_Y^V)$.
\end{proof}

\begin{remark}   \begin{enumerate}
\item Theorem \ref{t1} is valid for any dimension.
\item In general, the cone of sum forms will be a strict subset of the relative cone ${\mathcal C}_M^V$.  For example, consider $M=K3\#_{T^2_f}K3$, the fiber sum of two K3 surfaces along a fiber torus.  This has $b_2(M)=45$ which is also the dimension of the relative cone ${\mathcal C}_M^{T^2_f}$.  On the other hand, ${\mathcal C}_{K3}^{T^2_f}$ has dimension 22.  Hence the cone of sum forms must be a strict subset of the relative cone.  This indicates, that a precise study of the second homology of the symplectic sum $M=X\#_VY$ should be interesting, and we dedicate the rest of the section to this analysis.
\end{enumerate}\end{remark}

\subsubsection{The Second Homology  of $M=X\#_VY$}  We assume that $X,Y$ are 4-manifolds and that $V$ has trivial normal bundle.  The latter ensures that the
class of $V$ will exist in $H_2(X\#_V Y)$ after summing, albeit not
the particular copy of $V$ along which was summed. Denote the class
of $V$ by $f\in H_2(X\#_V Y)$.  In this section, we shall describe a "natural" basis of the second homology with respect to the fiber sum operation, which will allow us to efficiently construct and deconstruct cohomology classes on  $X\#_VY$.  

We begin by detailing the role of the maps involved in the symplectic sum in the structure of the second homology of $M=X\#_VY$.  The homology of $M$ can be analyzed by the Mayer-Vietoris sequences
for the triples $(X\sqcup_VY,X,Y)$ and $(X\#_VY,X\backslash
V,Y\backslash V)$ :
\begin{equation}
\begin{diagram}\label{rim}\dgARROWLENGTH=1em
\node{ H_2(V)}\arrow{e,t}{(j_X,j_Y)_*}\node{H_2(X)\oplus H_2(Y)}\arrow{e,t}{\phi}\node{H_2(X\sqcup_VY)}\\
\node{H_2(S_V)}\arrow{n}\arrow{e,t}{\lambda}\node{H_2(X\backslash V)\oplus H_2(Y\backslash V)}\arrow{n}\arrow{e,t}{\rho}\node{H_2(X\#_VY)}\arrow{n,l}{\psi}
\\
\node{R_V}\arrow{n}\node[2]{{\mathcal R}_{X\#_VY}}\arrow{n}
\end{diagram}
\end{equation}
The map $\lambda$ is induced by the canonical identification of the tubular neighborhoods and the normal bundles.  The map $\rho$ is identity on the classes which are supported away from $V$ and on classes supported near $V$ is defined by the gluing map $\varphi$, in particular by $\Theta$.  The map $\phi :H_2(X)\oplus H_2(Y)\rightarrow H_2(X\sqcup_VY)$ produces classes with the appropriate matching conditions on $V$ as determined by $\Theta$ in preparation for summing along $V$.  The map $\psi:H_2(X\#_VY)\rightarrow H_2(X\sqcup_VY)$ is induced by the gluing map $\varphi$, in particular by the embeddings $j_*$ and the isomorphism $\Theta$.  Then $\psi$ correctly decomposes classes in $X\#_VY$ in accordance with the symplectic gluing. The set ${\mathcal R}_{X\#_VY}$ is
completely determined by $R_V$ and an understanding of how these classes
bound in $M$. The outer columns are exact, for a detailed
discussion of the kernel $R_V$ see \cite{IP4}.

We are in the four dimensional setting, thus when we consider the Poincar\'e dual diagram to (\ref{rim}), we obtain in particular the following component:
\begin{equation}
\begin{diagram}\label{pd2}\dgARROWLENGTH=1em
\node{H^2(X\sqcup_VY)}\arrow{e,t}{\Psi}\arrow{s,l}{\Phi}\node{H^2(X\#_VY)}\arrow{e}\node{{\mathcal R}_{X\#_VY}^D}\arrow{e}\node{0}\\
\node{H^2(X)\oplus H^2(Y)}\\
\end{diagram}
\end{equation}
This is precisely the diagram used in the proof of Theorem \ref{t1}.  This motivates the  detailed discussion of the generators of the second homology which follows.

To explicitly describe the second homology of $M=X\#_VY$, we consider the following part of the Mayer-Vietoris sequences as before:
\begin{equation}
\begin{diagram}\label{star}\dgARROWLENGTH=1em
\node{H_2(X)\oplus H_2(Y)}\arrow{e,t}{\phi}\node{H_2(X\sqcup_VY)}\arrow{e,t}{\delta}\node{H_1(V)}\\
\node{H_2(X\backslash V)\oplus H_2(Y\backslash V)}\arrow{n}\arrow{e,t}{\rho}\node{H_2(M)}\arrow{n,l}{\psi}\arrow{e,t}{(\gamma,t)}\node{H_1(S_V)\simeq H_1(S^1)\oplus H_1(V)}\arrow{n,r}{\mu}\\
\node[2]{{\mathcal R}_{X\#_VY}}\arrow{n}\\
\end{diagram}
\end{equation}
Define the subgroups  $x=\rho(H_2(X\backslash V),0)$ and $y=\rho(0,H_2(Y\backslash V))$ of $H_2(M)$.  Elements in $x$ and $y$ are representable by submanifolds in $X$, $Y$ resp. which are supported away from the submanifold $V$.  Denote generators of $x$ and $y$ by $x_i$ and $y_i$.  Note that $x_i\cdot f=0=y_i\cdot f$ in the intersection form.

Define the subgroup $\Gamma=(\gamma,t)^{-1}(H_1(S^1),0)\simeq {\bf \mathbb Z}$.  Representatives $\gamma^M$ of this subgroup are submanifolds formed from submanifolds $\gamma^X\in X$ and $\gamma^Y\in Y$, these being supported in any neighborhood of $V$ and thus affected by the sum construction. The submanifolds $\gamma^*$ intersect $V$ nontrivially and $\gamma^*\cdot f=1$.  We denote the generator of this subgroup by $\gamma^M$.  Note that $\psi(\gamma^M)$ is always nontrivial: $\psi(\gamma^M)=(\gamma^X,\gamma^Y)$.

Define  $\tau=(\gamma,t)^{-1}(0,H_1(V))$ and ${\mathcal R}_M=\ker \psi$. The following holds:

\begin{lemma}  $\psi(\tau)= coker(\phi )$
\end{lemma}

\begin{proof}Let $g_i$ be generators of $H_1(V)$.  Then $\tau$ is generated by 
\[
\tau_i=(\gamma,t)^{-1}(0,g_i).
\]
The commutativity of (\ref{star}) shows
\begin{equation}
\delta\psi(\tau_i)=(\mu)(\gamma,t)\tau_i=\mu(0,t(\tau_i))=t(\tau_i)=g_i.
\end{equation}
Thus it follows from $g_i\ne 0$ that $\psi(\tau_i)\notin \ker \delta$.  Thus
\begin{equation}
\psi(\tau_i)\in H_2(X\sqcup_VY)/\ker \delta=H_2(X\sqcup_VY)/im \phi=coker(\phi).
\end{equation}
Hence $\psi(\tau)\subset coker(\phi)$.

From the definition of $coker (\phi)$ it follows that $coker (\phi)=H_2(X\sqcup_VY)/im \phi$ and hence any nontrivial element $c$ of the cokernel is supported in a neighborhood of $V$, but is not generated out of elements in $H_2(X)$ and $H_2(Y)$.  In particular, $c\cdot\gamma^*=0$.   Thus any lift $\tilde \tau$ of an element $c$ in the cokernel by $\psi$ to $H_2(M)$ has $\gamma(\tilde \tau)=0$.  Furthermore,
\begin{equation}
0\ne\delta(c)=\delta \psi(\psi^{-1}\tau)=(\mu)(\gamma,t)\tilde\tau=t(\tilde \tau),
\end{equation}
so $(\gamma,t)\tilde \tau\in 0\oplus H_1(V)$ is nontrivial.  Therefore, $\tilde \tau\in \tau$ and thus $coker (\phi)\subset \psi(\tau)$.
\end{proof}

Consider now the set ${\mathcal R}_{X\#_VY}$.  In particular, let us describe how objects in this set are generated from submanifolds of $X\backslash V$ and $Y\backslash V$. Define ${\mathcal R}_Y$ as the image of the map $i^Y_*\Delta: H_1(V)\rightarrow H_2(Y\backslash V)$ where $i^Y$ is the inclusion of $S_V$ into $Y\backslash V$ and $\Delta :H_1(V)\rightarrow H_2(S_V)$ stems from the Gysin sequence for the bundle $S_V\rightarrow V$.  If $V$ is 2 dimensional, then the map $\Delta$ is an injection.  Furthermore, consider for each simple closed curve $l$ in $V$ the preimage in $\partial N_YV$, this is a torus.  Such tori are called rim tori.  We restate a result in \cite{IP4}:

\begin{lemma}(Lemma 5.2,\cite{IP4}) Each element $R\in{\mathcal R}_Y$ can be represented by a rim torus.\end{lemma}

Under symplectic gluing, rim tori glue and are the elements of ${\mathcal R}_M$, in particular the elements $i_*^X\Delta l$ and $-i_*^Y\Delta l$ for some loop $l\in H_1(V)$ glue.  An example of this process is the generation of non-fiber tori in $K3$ when viewed as a sum $E(1)\#_{T^2_f}E(1)$ (See \cite{GS}).  These are then precisely the elements of ${\mathcal R}_{K3}=\{T^2_1,T^2_2\}$ and, in the same process,  ${\tau}=\{S^2_1,S^2_2\}$ is produced.  This accounts for the two new hyperbolic terms in the intersection form.  We observe the following

\begin{lemma}\label{K}Assume that $H_1(V)\rightarrow H_1(Y)$ is an injection and $V$ has trivial normal bundle.  Then $Y$ has no rim tori and $\tau=0={\mathcal R}_{X\#_VY}$.\end{lemma}

\begin{proof}To prove that $Y$ has no rim tori, it will suffice to show, that $i_*:H_2(S_V)\rightarrow H_2(Y\backslash V)$ is trivial on elements which are trivial under the map $\pi_*:H_2(S_V)\rightarrow H_2(V)$.  Therefore, consider the map $\xi: H_3(Y)\rightarrow H_2(S_V)$ where $S_V=V\times S^1$.  Let $W\in H_3(Y)$, then $\xi(W)=W\cap S_V=W\cap (V\times S^1)$.  In particular, $W\cap V\in H_1(V)$, thus by the injectivity assumption, if this intersection is non-trivial, it is non-trivial in $H_1(Y)$.  Therefore, the map $\xi$ maps $H_3(Y)$ onto the space generated by $\alpha\times S^1$ and $\beta \times S^1$, where $\alpha,\beta$ are generators of $H_1(V)$.  This space is the kernel of $\pi_*$ and the map $i_*$ is trivial on it.

Let us now consider ${\mathcal R}_{X\#_VY}$.  Elements of this set are constructed by symplectic gluing from elements in $X\backslash V$ and $Y\backslash V$, or equivalently, from classes in $H_2(X  \backslash V)$ and $H_2(Y\backslash V)$.  In particular, considering \ref{rim}, only classes in $R_V$ are relevant, these are precisely those classes which do not map trivially to $H_2(V)$.  As we have seen above, our assumption implies that $H_2(S_V)\rightarrow H_2(Y\backslash V)$ is trivial.  Hence every element in ${\mathcal R}_{X\#_VY}$ would be trivial.

Thus ${\mathcal R}_{X\#_VY}=0$.  Furthermore,  $\tau=0$ is clear from the assumption.
\end{proof}

The previous discussion allows us to explicitly state a set of generators for $H_2(M)$:

\begin{itemize}
\item $\{f\}$ is the fiber class present in both $X$ and $Y$, in our case this is the class of $V$;
\item $\{x_i\}$, $\{y_i\}$;
\item $\{k_i\}\subset {\mathcal R}_{X\#_VY}$, generators which are represented by submanifolds mapping to 0 in $X\sqcup_VY$;
\item $\{\gamma\}$ generated out of elements of the homology of both $X$ and $Y$, e.g $[\Sigma_2]$ from copies of $T^2$ in Example \ref{T4T4}.  Note that this is the origin for the non-surjectivity of the map $\psi$:  $\psi[\Sigma_2]$ will always have a fixed relative orientation of the two copies of $T^2$ into which $\Sigma_2$ degenerates.  Thus the pairing of the tori with opposite orientation will not lie in the image of $\psi$.
\item $\{\tau_i\}\subset \tau$;  these objects will persist in $X\sqcup_VY$ and hence contribute to its homology as well.
\end{itemize}

Given this set of generators, we can explicitly state how an element in the cone of sum forms decomposes:  Given $\alpha=\sum_i a_iX_i+b_iY_i+cF+g\Gamma+e_i{\mathcal R}_i+t_i T_i\in {\mathcal C}^{sum}$ and taking the Poincar\'e dual basis  of the one given above, we obtain two forms $\alpha_X=\sum_i a_iX_i+c^XF^X+g\Gamma^X$ and $\alpha_Y=\sum_i b_iY_i+c^YF^Y+g\Gamma^Y$ on $X$ resp. $Y$.  Note that this is ultimately a direct result of Theorem \ref{t1}.

\subsubsection{Good sums}

If we know the relative cones of $X$ and $Y$, then,  considering \ref{rim}, we should obtain information on the structure of the relative cone on $M=X\#_VY$ by using the Poincar\'e duals of the maps $\phi$ and $\psi$.  For this to work nicely, one needs $\phi$ to be surjective and $\psi$ to be injective.  We thus make the following definition:

\begin{definition}A symplectic sum $M=X\#_VY$ is called \underline{good} if $\phi$ is surjective and $\psi$ is injective.\end{definition}

This statement is equivalent to $\tau=0={\mathcal R}$, and Lemma \ref{K} provides a simple criterion to check this.

\begin{theorem}\label{t2} Suppose $M=X\#_VY$ is good and  $V$ has trivial normal bundle.   If   for $X,Y$,
\begin{equation}
{\mathcal C}^V_*=\{\alpha\in {\mathcal P}_*\;\vert\;\alpha\cdot [V]^D> 0\},
\end{equation}
then
\begin{equation}
{\mathcal C}^{sum}_{X\#_VY}= \{\alpha\in {\mathcal P}_M\;\vert\;\alpha\cdot [V]^D>
0\}.
\end{equation}
Consequently ${\mathcal C}_M^V=\{\alpha\in {\mathcal P}_M\;\vert\;\alpha\cdot [V]^D>
0\}$.

\end{theorem}

\begin{proof}The second result is immediate:  Theorem \ref{t1} and Lemma \ref{incl} show that ${\mathcal C}_{X\#_VY}^{sum}\subset {\mathcal C}_M^V\subset \{\alpha\in {\mathcal P}_M\;\vert\;\alpha\cdot [V]^D>0\}$.

The first result follows, if we can show $ \{\alpha\in {\mathcal P}_M\;\vert\;\alpha\cdot [V]^D>0\}\subset {\mathcal C}_{X\#_VY}^{sum}$.

We proceed as remarked above, using Theorem \ref{t1}:  Taking the Poincar\'e dual basis  of the one given above, we can write each $\alpha\in \{\alpha\in {\mathcal P}_M\;\vert\;\alpha\cdot [V]^D>0\}$ as
\begin{equation}
\label{dec}\alpha=\sum_i a_iX_i+b_iY_i+cF+g\Gamma+e_i{\mathcal R}_i+t_i T_i;\;\;\;g>0.
\end{equation}
As ${\mathcal R}=0=\tau$, the last two terms drop.

We must now show, that $\alpha\in {\mathcal C}_{X\#_VY}^{sum}$.  We thus choose a possible pair of classes $\alpha_X$ and $\alpha_Y$ in $H^2(X)$ resp. $H^2(Y)$ as determined by Theorem \ref{t1} and show that this can be done in such a way as to ensure that they are representable by a relative symplectic form. We first determine a relation which preserves the volume.  In the following, we show how to choose this pair, so that they lie in their respective relative cones ${\mathcal C}_*^V$.  Hence the class $\alpha$ can be obtained by summing two classes representable by relative symplectic forms and  thus, by Gompf's result, $\alpha\in {\mathcal C}_{X\#_VY}^{sum}$.

Choose the candidates for classes summing to  $\alpha$ as follows:
\begin{equation}
\begin{array}{ll}
\alpha_X&=\sum_i a_iX_i+c^XF^X+g\Gamma^X\in H^2(X)\\
\alpha_Y&=\sum_i b_iY_i+c^YF^Y+g\Gamma^Y\in H^2(Y)
\end{array}
\end{equation}
where $F^*$ and $\Gamma^*$ are the Poincar\'e duals on $X,Y$.  The coefficient $g$ must be the same for both, as $g=\alpha([V])=\alpha_X([V])=\alpha_Y([V])$.  The class $F$ has $F^2=0$ due to the triviality of the normal bundle of $V$, similarly $(F^*)^2=0$.  The volume of each of these is
\begin{equation}
\alpha^2= (\sum a_iX_i)^2+(\sum b_iY_i)^2+( g\Gamma)^2+\end{equation}\[+2\sum\left(a_iX_ig\Gamma+b_iY_ig\Gamma\right)+cFg\Gamma\]
and
\begin{equation}
\alpha_X^2=(\sum a_iX_i)^2+ g^2(\Gamma^X)^2+2\sum\left(a_iX_ig\Gamma^X\right)+c^XF^Xg\Gamma^X.
\end{equation}
Thus the difference of the volumes is calculated to be
\begin{align}
\alpha^2-\alpha_X^2-\alpha_Y^2
&\label{e1}=\left(g\Gamma\right)^2-\left( g\Gamma^X\right)^2-\left( g\Gamma^Y\right)^2\\
&\label{e2}+2\left(\sum a_iX_ig\Gamma -\sum a_iX_ig\Gamma^X\right)\\
&\label{e3}+2\left(\sum b_iY_ig\Gamma-\sum b_iY_ig\Gamma^Y\right)\\
&\label{e4}+2\left( cFg\Gamma-c^XF^Xg\Gamma^X-c^YF^Yg\Gamma^Y\right)
\end{align}
Note the following:  The morphism $\Psi: H^2(X\sqcup_VY)\rightarrow H^2(M)$ relates the intersection forms, giving the following relations:
\begin{enumerate}
\item $(\Gamma^X)^2 +(\Gamma^Y)^2=(\Gamma^X\oplus \Gamma^Y)^2=\Psi((\Gamma^X\oplus \Gamma^Y)^2)=\Psi(\Gamma^X\oplus \Gamma^Y)^2=\Gamma^2$
\item $\alpha_iX_i\beta\Gamma=\Psi(\alpha_iX_i(\beta\Gamma^X\oplus \beta\Gamma^Y))=\alpha_iX_i\beta\Gamma^X$
\end{enumerate}
Applying these relations, it follows immediately that \ref{e1} is trivial,
\begin{equation}
\ref{e2}\Rightarrow a_iX_ig\Gamma- a_iX_ig\Gamma^X =a_i(g-g)X_i\Gamma^X=0
\end{equation}
and analogously for \ref{e3} and $Y$.  Furthermore, (\ref{e4}) becomes
\begin{equation}
cFg\Gamma-c^XF^Xg\Gamma^X-c^YF^Yg\Gamma^Y\end{equation}\[=cFg\Gamma-\Psi(c^XF^Xg\Gamma^X+c^YF^Yg\Gamma^Y)\]\[=cFg\Gamma-\Psi(c^XF^X+c^YF^Y)g\Gamma.
\]
The condition for this to vanish is
\begin{equation}
\Psi(c^XF^X+c^YF^Y)=cF.\label{dcond}
\end{equation}
Thus, by choosing  $c=c^X+c^Y$  we preserve the volume.

Now, we must show that this can be done in such a way, as to ensure $\alpha_*^2>0$. Choose $c^X,c^Y$ so that volume is preserved.  Then $\alpha^2=\alpha^2_X+\alpha_Y^2>0$, and we may assume  $\alpha_X^2>0$.  This holds true for any choice of $c^*$ satisfying \ref{dcond}.

Squaring $\alpha_X$ and denoting $B=\sum a_iX_i+g \Gamma^X$, we obtain
\begin{equation}
\label{par}f(c^X)=\alpha_X^2=B^2 +2B\cdot F^X c^X+ (c^X)^2(F^X)^2=B^2 +2B\cdot F^X c^X.
\end{equation}
We can always solve $f(c^X)=\rho$ for any $\rho>0$.  Thus we can ensure that $\alpha^2>\alpha_X^2>0$ holds.  Then also $\alpha_Y^2=\alpha^2-\alpha_X^2>0$ holds, and thus each $\alpha_*$  must lie in ${\mathcal C}_*^V$, hence  $\tilde\alpha=(\alpha_X,\alpha_Y)\in {\mathcal  C}_{X\sqcup_VY}$ by definition of this set.  Thus $\{\alpha\in {\mathcal P}_*\;\vert\;\alpha\cdot [V]^D> 0\}\subset {\mathcal C}_{X\#_VY}^{sum}$.
\end{proof}

This result is of particular interest, as it shows that good sums preserve the structure of the relative cone.  Thus, if $X,Y$ have relative cones as assumed in the Theorem and the sum is good, we can apply this result repeatedly to obtain the relative cone of $nX\#_VmY$:
\begin{equation}
\label{bigsum}{\mathcal C}_{nX\#_VmY}^V=\{\alpha\in {\mathcal P}_{nX\#_VmY}\;\vert\;\alpha\cdot [V]^D> 0\}.
\end{equation}

\section{\label{examples}Symplectic cone of certain $T^2-$fibrations}

\subsection{$T^2\times \Sigma_g$}

We now show that Theorem \ref{t2} can be applied to $T^2\times
\Sigma_g$.  The results of the previous section assume two things: a
certain form of the relative symplectic cone and that the sum be
good.

Fix $Y=T^2\times \Sigma_k$. Thus by Lemma \ref{K} we don't need to
verify the condition ${\mathcal R}=0=\tau$
when applying Theorem \ref{t2}, i.e. all sums $X\#_VY$ are good.

The following result follows immediately:

\begin{theorem} \label{SnT4}Let  $M=T^2\times \Sigma_k$. Then
\begin{equation}
{\mathcal C}_M^{T^2_f}=\{\alpha\in {\mathcal P}_M\;\vert\;\;\alpha\cdot [T^2_f]^D> 0\}.
\end{equation}
Consequently,  ${\mathcal C}_M={\mathcal P}^{c_1(M)}\cup {\mathcal
P}^{-c_1(M)}$.
\end{theorem}

\begin{proof}We proceed by induction:  Let $M=T^4$.  Then the result holds
due to Lemma \ref{Torus} and Lemma \ref{VD}. Summing repeatedly
we obtain $$M=T^2\times \Sigma_k= T^4\#_{T^2_f}(T^2\times
\Sigma_{k-1}).$$ Using the induction hypothesis, which ensures that
$${\mathcal C}_{T^2\times \Sigma_{k-1}}^{T^2_f} =\{\alpha\in {\mathcal
P}_{T^2\times \Sigma_{k-1}}\;\vert\;\;\alpha\cdot [T^2_f]^D> 0\},$$
the result now follows from Theorem \ref{t2} (see \ref{bigsum}) and
Lemma \ref{VD}. 
\end{proof}

This result also follows from results in \cite{FV1} and \cite{FV2}.

\begin{remark} (Fibered symplectic forms) Every class in ${\mathcal C}_M$
can be represented by a symplectic form which restricts to a
symplectic form on the fibers of $M$. Denote the set of such forms
by ${\mathcal S}$.
 Then this set is contractible (and nonempty) (Thm. 1.4, \cite{Go5}).

 See also McDuff \cite{M1}.
\end{remark}

 \subsection{$X\# (T^2\times \Sigma_k)$}

 In the following we allow the fibration to
have  singular or multiply covered fibers. If we sum along a generic
fiber, avoiding these special fibers, we find no obstruction to
applying the methods developed above.

\begin{theorem}Let $X$ be a minimal elliptic K\"ahler surface with $p_g=0$ and $M=X\#_{T^2_f}Y$.  Then Conjecture \ref{conj} holds, i.e.$$C_M={\mathcal P}^{c_1(M)}\cup {\mathcal P}^{-c_1(M)}.$$

\end{theorem}

\begin{proof}This reduces to an application of Theorem \ref{t2} and Lemmas
\ref{K} and \ref{VD} together with Theorem \ref{li1} to
show the result.  
\end{proof}

{\bf Remark:} The manifold $X$ could be an Enriques surface, a hyperelliptic surface or a Dolgachev surface.  It should also be possible to use this method to determine the symplectic cone of $X\#_{T^2_f}(T^2\times \Sigma_g)$ with $X$ a minimal symplectic manifold with $b^+=1$ and which contains a square 0 torus $T^2_f$.

\section{Appendix}

Let $V$ be a fixed smooth codimension 2 submanifold of a symplectic manifold $(X,\omega)$. Let $\mathcal J_V$ be the set of almost complex structures compatible with $\omega$ such that $V$ is pseudoholomorphic for each $j\in\mathcal J_V$.  We wish to show that $\mathcal J_V$ has a rich enough structure to allow for genericity statements for $J$-holomorphic curves.  These results are presumably known to experts in the field, they use methods in \cite{McS} and \cite{U}, we include them for completeness.  Let $A\in H_2(X)$ be any class, except that in the case $A^2=0=K_\omega(A)$ the class $A$ should be indivisible.  We begin by defining a universal space which we shall use throughout this section:  Fix a closed compact Riemann surface $\Sigma$.  The universal model $\mathcal U$ is defined as follows: This space will consist of Diff$(\Sigma)$ orbits of a 4-tuple $(i,u,J,\Omega)$ with
\begin{enumerate}
\item $u:\Sigma\rightarrow X$ an embedding off a finite set of points from a Riemann surface $\Sigma$ such that $u_*[\Sigma]=A$ and $u\in W^{k,p}(\Sigma,X)$ with $kp>2$,
\item $\Omega\subset X$ a set of $m$ distinct points (with $\Omega=\emptyset$ if $m\le 0$) such that $\Omega\subset u(\Sigma)$,
\item $i$ a complex structure on $\Sigma$ and $J\in \mathcal J_V$.
\end{enumerate}
Note that every map $u$ is locally injective. 

In order to show the necessary genericity results, we will call upon the Sard-Smale Theorem.  This will involve the following technical difficulty:  The spaces $\mathcal J_V$ and any subsets thereof which we will consider are not Banach manifolds in the $C^\infty$-topology.  However, the results we wish to obtain are for smooth almost complex structures.  In order to prove our results, we need to apply Taubes trick (see \cite{T} or \cite{McS}), which replaces the smooth spaces by $C^l$-almost complex structures and apply the Sard-Smale Theorem in that setting.  Then, one constructs a countable collection of  sets, whose intersection is the generic set of smooth structures, and shows that each of these is dense and open by explicit argumentation in the space of smooth structures.  We will not go through this technical step but implicitly assume this throughout the section, details can be found in Ch. 3 of \cite{McS}.

\begin{lemma}\label{genericA}
Let $A\in H_2(X,\mathbb Z)$, $A\ne [V]$, and $k(A)=\frac{1}{2}(A^2-K_\omega(A))\ge 0$.  Let $\Omega$ denote a set of $k(A)$ distinct points in $X$.  Denote the set of pairs $(J,\Omega)\in \mathcal J_V\otimes X^{k(A)}$ by $\mathcal I$.  Let $\mathcal J_V^A$ be the subset of pairs $(J,\Omega)$ which are non-degenerate for the class $A$ in the sense of Taubes (\cite{T}). Then $\mathcal J_V^A$ is a set of second category in $\mathcal I$.
\end{lemma}

The term non-degenerate as defined by Taubes states that for the pair $(J,\Omega)$ the linearization of the operator $\overline\partial$ for any $J$-holomorphic submanifold of $X$ containing the set $\Omega$ has trivial cokernel.  Note that the universal model excludes multiple covers of the submanifold $V$ in the case that $A=a[V]$ for $a\ge 2$, and we can thus assume that any map $u:\Sigma\rightarrow M$ satisfies $u(\Sigma)\not\subset V$.

\begin{proof}To prove this statement, we will define a map $\mathcal F$ from a universal model $\mathcal U$ to a bundle with fiber $W^{k-1,p}(\Lambda^{0,1}T^*\Sigma\otimes u^*TX)$ and show that it is submersive at its zeroes.  Then we can apply the Sard-Smale theorem to obtain that $\mathcal J_V^A$ is of second category.

Define the map $\mathcal F$ as $(i,u,J,\Omega)\mapsto \overline\partial_{i,J}u$.  Then the linearization at a zero $(i,u,J,\Omega)$ is given as
\begin{equation}
{\mathcal F}_*(\alpha,\xi,Y)=D_u\xi + \frac{1}{2}(Y\circ du\circ i + J\circ du\circ \alpha)
\end{equation}
where $D_u$ is Fredholm, $Y$ and $\alpha$ are variations of the respective almost complex structures.  

Consider $u\in \mathcal U$ such that there exists a point $x_0\in \Sigma$ with $u(x_0)\in X\backslash V$ and $du(x_0)\ne 0$ (The second condition is satisfied almost everywhere, as $u$ is a $J$-holomorphic map.).  Then there exists a neighborhood $N$ of $x_0$ in $\Sigma$ such that 
\begin{enumerate}
\item $du(x)\ne 0$, 
\item $u(x)\not\in V$ for all $x\in N$.
\end{enumerate}
In particular, we know that the map $u$ is locally injective on $N$.  Furthermore, we can find a neighborhood in $N$, such that there are no constraints on the almost complex structure $J\in \mathcal J_V$, i.e. this neighborhood does not intersect $V$.  In particular, $Y$ can be chosen as from the set of $\omega$-tame almost complex structures.  Denote this open set by $N$ as well.

Let $\eta\in $coker$ \mathcal F_*$.  Consider any $x\in N$ with $\eta(x)\ne 0$.  Then Lemma 3.2.2, \cite{McS}, provides a matrix $Y_0$ with the properties
\begin{itemize}
\item $Y_0=Y_0^T=J_0Y_0J_0$ with $J_0$ the standard almost complex structure in a local chart and 
\item $Y_0[ du(x)\circ i(x)]=\eta(x)$.
\end{itemize}
On $N$ choose any variation $Y$ of $J$ such that $Y(u(x))=Y_0$.  Then define the map $f:N\rightarrow \mathbb R$ by $\langle Y\circ du\circ i,\eta\rangle$.  Note that $f(x)>0$ by definition of $Y$.  Therefore, we can find an open set $N_1$ in $N$ such that $f>0$ on that open set. 
Using the local injectivity of the map $u$ and arguing as in \cite{McS}, we can find a neighborhood $N_2\subset N_1$ and a neighborhood $U\subset M$ of $u(x_0)$ such that $u^{-1}(U)\subset N_2$.  Choose a cutoff function $\beta$ supported in $U$ such that $\beta(u(x))=1$.  Hence in particular
\begin{equation}
\int_\Sigma \langle \mathcal F_*(0,0,\beta Y),\eta\rangle>0
\end{equation}  
and therefore $\eta(x)=0$.  This result holds for any $x\in N$, therefore $\eta$ vanishes on an open set.

As we have assumed $\eta\in $coker$ \mathcal F_*$, it follows that
\[
0\;=\;\int_\Sigma \langle \mathcal F_*(0,\xi,0),\eta\rangle\;=\;\int_\Sigma\langle D_u\xi,\eta\rangle
\]
for any $\xi$.  Then it follows that $D_u^*\eta =0$ and $0=\triangle \eta + l.o.t.$.  Therefore Aronszajn's theorem allows us to conclude that $\eta=0$ and hence $\mathcal F_*$ is surjective.

Thus we have the needed surjectivity for all maps admitting $x_0$ as described above: $u(x_0)\not\in V$ and $du(x_0)\ne 0$.  As stated before, this last condition is fulfilled off a finite set of points on $\Sigma$.  The first holds for any map $u$ in class $A$ as we have assumed that $A\ne [V]$.

Now apply the Sard-Smale theorem to the projection onto the last two factors of $(i,u,J,\Omega)$.
\end{proof}

Define the set ${\mathcal K}^A_V(J,\Omega)$ to be the set of $J$-holomorphic submanifolds which are abstractly diffeomorphic to a Riemann surface $\Sigma$, contain the set $\Omega$, represent the class $A$ and $J\in\mathcal J_V$.  Then the same methods as in the above proof  together with index calculations of the projection operator onto the last two factors lead to the following results:  If $m>k(A)$ or $m<0$, then ${\mathcal K}^A_V(J,\Omega)$ is empty for generic $(J,\Omega)$, if $m=k(A)$, then ${\mathcal K}^A_V(J,\Omega)$ is a smooth 0-dimensional manifold for generic $(J,\Omega)$.  In particular, there exists a set of second category in $\mathcal J_V$, such that if $k(A)\ge 0$ then any pseudoholomorphic submanifold in class $A$ meets a generic set of at most $k(A)$ distinct points.

As we have seen in the above proof, for the class $A=[V]$ which may have representatives which do not lie outside of $V$, we must be careful.  In particular, it is conceivable, that the particular submanifold $V$ chosen may not be generic in the sense of Taubes, i.e. the set $\mathcal J_V$ may contain almost complex structures for which the linearization of $\overline \partial_J$ at the embedding of $V$ is not surjective.  The rest of this section addresses this issue.  We begin by showing that the cokernel of the linearization of the operator $\overline\partial_J$ at a $J$-holomorphic embedding of $V$ has the expected dimension:

Let $j$ be an almost complex structure on $V$.  Define $\mathcal J_V^j=\{J\in\mathcal J_V\;\vert \;J\vert_V=j\}$ and call any $J$-holomorphic embedding of $V$ for $J\in  \mathcal J_V^j$ a $j$-holomorphic embedding.

\begin{lemma}\label{genericV}Fix a $j$-holomorphic embedding $u:(\Sigma,i)\rightarrow (X,J)$ for some $J\in\mathcal J_V^j$.    If $k([V])\ge 0$, then there exists a set $\mathcal J_V^{g,j}$ of second category  in $\mathcal J_V^j$ such that for any $J\in \mathcal J_V^{g,j}$ the linearization of $\overline \partial_{i,J}$ at the embedding $u$ is surjective.  If $k([V])<0$, then then there exists a set $\mathcal J_V^{g,j}$ of second category  in $\mathcal J_V^j$ such that the submanifold $V$ is rigid in $X$. 
\end{lemma}

Let us explain the structure of the proof before giving the exact proof.  We follow ideas of section 4, \cite{U}.  We need to show that for a fixed embedding $u:\Sigma\rightarrow X$ of $V$ the linearization $\mathcal F_*$ of $\overline \partial_{i,J}$ at u has a cokernel of the correct dimension for generic $J\in \mathcal J_V^j$.  To do so, we will consider the operator $\mathcal G(\xi,\alpha,J):= \mathcal F_*(\alpha,\xi,0)$ at $(i,u,J,\Omega)$.  We will show that the kernel of the linearization $\mathcal F_*$ for non-zero $\xi$ has the expected dimension for generic $J$ and hence the linearization of $\overline \partial_{i,J}$ at $u$ also has the expected dimension.  Note also, that for any $J\in\mathcal J_V^j$, the map $u$ is $J$-holomorphic.

What is really going on in this construction?  The operator $\mathcal F$ is a section of a bundle over $\mathcal U$, as described above.  Further, we consider a map $\mathcal U\rightarrow \mathcal J_V^j$.  In this map, we fix a "constant section" $(u,j)$, i.e. we consider the structure of the tangent spaces along a fixed map $u$ where we do not let the almost complex structure along $V$ vary.  On he other hand, it is only this structure $j$ which makes $u$ pseudoholomorphic.  Hence fixing $(u,j)$ is akin to considering a constant section in the bundle $\mathcal U\rightarrow \mathcal J_V^j$.  In particular, we are only interested in the component of the tangent space along this section  which corresponds to the tangent space along the moduli space $\mathcal M=\mathcal F^{-1}(0)$, as this will give us insight into the dimension of $\mathcal M$.  Along $(u,j)$, this is precisely the component of the kernel of $\mathcal F_*$ with $Y=0$ as the complex structure is fixed on $V$, i.e. the set of pairs $(\xi,\alpha)$ such that $\mathcal F_*(\xi,\alpha,0)=0$, which corresponds to exactly the zeroes of $\mathcal G$.  When considering the zeroes of the map $\mathcal G$ viewed over $\mathcal J_V^j$, we find that this is a collection of finite dimensional vector spaces.  We may remove any part of these spaces, so long as we leave an open set, which is enough to allow us to determine the dimension of the underlying vector spaces.  Hence, removing $\xi=0$, a component along which we cannot use our methods to determine the dimension of the kernel, still leaves a large enough set to be able to determine the dimension of the moduli space $\mathcal M$.  We therefore want to show that the kernel of the linearization $\mathcal F_*$ for non-zero $\xi$ which is the zero set of $\mathcal G$ for non-zero $\xi$ has the expected dimension $\max\{k([V]), 0\}$ for generic $J$.

\begin{proof}The operator $\mathcal G$ is defined as 
\[
W^{1,p}(u^*TX)\times H^{0,1}_i(T_{\mathbb C}\Sigma)\times \mathcal J_V^j\rightarrow L^p(u^*TX\otimes T^{0,1}\Sigma)
\]
\[
(\xi,\alpha,J)\mapsto D_u^J\xi + \frac{1}{2}J\circ du\circ \alpha
\]  
where the term $D_u^J=\frac{1}{2}(\nabla \xi+ J\nabla\xi\circ i)$ for some $J$-hermitian connection $\nabla$ on $X$, say for example the Levi-Civita connection associated to $J$.  Note that we could define this operator also for a smooth embedding $u:\Sigma\rightarrow X$, but that we have fixed the almost complex structure on $V$ and therefore $\overline\partial_{i,J}u=0$ for any $J\in \mathcal J_V^j$.

Let $(\xi,\alpha,J)$ be a zero of $\mathcal G$.  Linearize $\mathcal G$ at $(\xi,\alpha,J)$:
\[
\mathcal G_*(\gamma,\mu,Y)=D_u^J\gamma +\frac{1}{2}\nabla_\xi Y\circ du\circ i+\frac{1}{2}J\circ du\circ \mu.
\]
As stated above, we assume nonvanishing $\xi$, hence we can assume that $\xi\ne 0$ on any open subset.  Let $\eta\in $coker $ \mathcal G_*$.  Let $x_0\in \Sigma$ be a point with $\eta(x_0)\ne 0\ne \xi(x_0)$.  In a neighborhood of $u(x_0)\in V$ the tangent bundle $TX$ splits as $TX=N_V\oplus TV$ with $N_V$ the normal bundle to $V$ in $X$.  With respect to this splitting, the map $Y$ has the form
\[
y=\left(\begin{array}{cc} a&b\\0&0\end{array}\right)
\]
with all entries $J$-antilinear and $b\vert_V=0$, thus ensuring that $V$ is pseudoholomorphic and accounting for the fact that we have fixed the almost complex structure along $V$.  Thus $\nabla_\xi Y$ can have a similar form, but with no restrictions on the vanishing of components along $V$.  In particular, assuming $\eta$ projected to $N_V$ is non-vanishing, we can choose 
\[
\nabla_\xi Y=\left(\begin{array}{cc} 0&B\\0&0\end{array}\right)
\]
at $x_0$ such that $B(x_0)[du(x_0)\circ i(x_0)](v)=\eta^{N_V}(x_0)(v)$ and $B(x_0)[du(x_0)\circ i(x_0)](\overline v)=\eta^{N_V}(x_0)(\overline v)$ for a generator $v\in T^{1,0}_{x_0}\Sigma$ and where $\eta^{N_V}$ is the projection of $\eta$ to  $N_V$.  Then, using the same universal model as in the previous Lemma, we can choose neighborhoods of $x_0$ and a cutoff function $\beta$ such that 
\[
\int_\Sigma \langle \mathcal G_*(0,0,\beta Y),\eta\rangle >0
\]
and thus any element of the cokernel of $\mathcal G_*$ must have $\eta^{N_V}=0$.  An argument in \cite{U} shows that the projection of $\eta$ to $TV$ must also vanish.  Therefore the map $\mathcal G_*$ is surjective at the embedding $u:\Sigma\rightarrow V$.  

Thus the set $\{(\xi,\alpha,J)\vert \mathcal G(\xi,\alpha,J)=0,\;J\in\mathcal J_V^j,\;\xi\ne 0\}$ is a smooth manifold and we may project onto the last factor.  Then applying Sard-Smale, we obtain a set $\mathcal J_V^{g,j}$ of second category in $\mathcal J_V^j$, such that for any $J\in \mathcal J_V^{g,j}$, the kernel of the linearization of $\overline\partial $ at non-zero perturbations $\xi$ of the map $u$ is a smooth manifold of the expected dimension.  In the case $k([V])\ge 0$, this however implies that $\mathcal F_*$ at $(i,u,J,\Omega)$ is surjective.  Therefore, we have found a set $\mathcal J_V^{g,j}$ of second category in $\mathcal J_V^j$ such that the linearization of $\overline \partial_{i,J}$ at u is surjective at all elements of $\mathcal J_V^{g,j}$.  

If however $k([V])<0$, then this kernel is generically empty.  This implies the rigidity of the embedding $u$ of $V$.

\end{proof}

We have thus shown, that for a fixed embedding we can find a generic set of almost complex structures among those making the embedding pseudoholomorphic, such that the linearization of $\overline\partial$ at u has cokernel of the expected dimension.   We would like to state a similar result for the space $\mathcal K_V^{[V]}(J,\Omega)$ as we stated for the class $A$.  In order to do so, note that in our model, we consider orbits under the action of Diff$(\Sigma)$.  Hence, given any two embeddings $u:(\Sigma,i)\rightarrow (X,J)$ and $v:(\Sigma,i)\rightarrow (X,\tilde J)$  of $V$ for $J,\tilde J\in\mathcal J_V^j$, there exists a $\phi\in$ Diff$(\Sigma)$ such that $u=v\circ \phi$, i.e. $u$ and $v$ correspond to the same point in the universal model $\mathcal U$.

For every almost complex structure $j$ on $V$ the previous results provide the following:
\begin{enumerate}
\item A set $\mathcal J_V^{g,j}$ of second category in $\mathcal J_V^j$ with the property that the linearization of the operator $\overline\partial$ at a fixed $j$-holomorphic embedding of $V$  is surjective ($k([V])\ge 0$) or is injective ($k([V])<0$).
\item Up to  a map $\phi\in$ Diff$(\Sigma)$,  there is a unique $j$-holomorphic embedding of $V$ for all $J\in \mathcal J_V^j$.
\end{enumerate}

Therefore, consider the following set:
\[
\mathcal J_V^g=\bigcup_{j} \mathcal J_V^{g,j}\subset \bigcup_{j} \mathcal J_V^j=\mathcal J_V.
\]

Note that we $\mathcal J_V^g$ is actually a disjoint union of sets.  The following properties hold:  
\begin{enumerate}
\item The set $\mathcal J_V^g$ is dense in $\mathcal J_V$.
\item The linearization of the operator $\overline\partial$ at a fixed $j$-holomorphic embedding of $V$  is surjective ($k([V])\ge 0$) or is injective ($k([V])<0$) for any $J\in \mathcal J_V^g$.
\item Up to  a map $\phi\in$ Diff$(\Sigma)$,  there is a unique $j$-holomorphic embedding of $V$.
\end{enumerate}
We can now state the final result concerning genericity that we will need:

\begin{lemma}\label{jV}  Let $\Omega$ denote a set of $k([V])$ distinct points in $X$.
\begin{enumerate}
\item $k([V])\ge 0$: Denote the set of pairs $(J,\Omega)\in \mathcal J_V\times X^{k([V])}$ by $\mathcal I$ (with $\Omega=\emptyset$ if $k([V])\le 0$).  Let $\mathcal J_{[V]}$ be the subset of pairs $(J,\Omega)$ which are non-degenerate for the class $A$ in the sense of Taubes (\cite{T}).  Then $\mathcal J_{[V]}$ is dense in $\mathcal I$.
\item $k([V])<0$: There exists a dense set $\mathcal J_{[V]}\subset \mathcal J_V$ such that $V$ is rigid, i.e. there exist no pseudoholomorphic deformations of $V$ and there are no other pseudoholomorphic maps in class $[V]$. 
\end{enumerate}
\end{lemma}

\begin{proof}To begin, we will replace the set $\mathcal J_V\times X^{k([V])}$ by $\mathcal J^g_V\times X^{k([V])}$ which is a dense subset, as seen from the previous remarks.  Further, for any $(J,\Omega)\in \mathcal J^g_V\times X^{k([V])}$, we have surjectivity or injectivity of the linearization at the embedding of $V$.

Consider the case $k([V])\ge 0$.  Fix a $j$ on $V$.  Then consider the set $\mathcal J_V^{g,j}$ provided by Lemma \ref{genericV}.  The linearization at the embedding of $V$ is surjective for any $J\in \mathcal J_V^{g,j}$.  For any element $(i,u,J,\Omega)$ of $\mathcal U$ with $u(\Sigma)\not\subset V$ representing the class $[V]$ and $J\in \mathcal J_V^{g,j}$, arguments as in the proof of Lemma \ref{genericA} provide the necessary surjectivity.  Therefore, there exists a further set $\mathcal J_{[V]}^{g,j}$ of second category in $\mathcal J_V^{g,j}\times X^{k([V])}$ such that any pair $(J,\Omega)\in\mathcal J_{[V]}^{g,j}$ is nondegenerate.  

Define $\mathcal J_{[V]}=\bigcup_{j}\mathcal J_{[V]}^{g,j}$.  This is a dense subset of $\mathcal J^g_V\times X^{k([V])}$ such that any pair $(J,\Omega)\in\mathcal J_{[V]}$ is nondegenerate. 

If $k([V])<0$, then restrict to $ \mathcal J^g_V$ as well.  Thereby we have already ensured that $V$ is rigid.  Now apply the proof of Lemma \ref{genericA} to the universal model $\mathcal U$, which we modify to allow only maps $u:(\Sigma,i)\rightarrow (X,J)$ such that $u(\Sigma)\not\subset V$.  Then we can find a set $\mathcal J_{[V]}$ of second category in  $ \mathcal J^g_V$ such that there exist no maps in class $[V]$ other than the embedding of $V$.
\end{proof}

Note that by results of Taubes, if $k([V])<0$, the set $ \mathcal J_{[V]}$ is a set of first category in $\mathcal J$.  Further, it is not clear, whether it is possible to improve the denseness statement to include openness in $\mathcal J_V$.  In the case $k([V])\ge 0$, this is also unclear.

Furthermore, if $k([V])\ge 0$, then we have shown that the set $\mathcal K_V^{[V]}(J,\Omega)$ has the desired properties, i.e. for a dense set of pairs $(J,\Omega)$ $\mathcal K_V^{[V]}(J,\Omega)$ is a smooth 0-dimensional manifold unless $m>k([V])$, in which case it is generically empty.

Similar results have been proven by Jabuka in \cite{J}.  However, that result only provides an isotopic copy of $V$ in the case $k([V])\ge 0$. 



\begin{thebibliography}{99}

\bibitem{BPV}Barth, Wolf P. ;  Hulek, Klaus ;  Peters, Chris A. M. ;  Van de Ven, Antonius . {\it Compact complex surfaces.
Second edition.}
Ergebnisse der Mathematik und ihrer Grenzgebiete. 3. Folge. A Series of
 Modern Surveys in Mathematics [Results in Mathematics and Related Areas. 3rd
 Series. A Series of Modern Surveys in Mathematics], 4. Springer-Verlag, Berlin,  2004. xii+436 pp.  
		

\bibitem{Bi} Biran, Paul . {\it A stability property of symplectic packing.}
 Invent. Math.  136  (1999),  no. 1, 123--155.
\bibitem{BL}Bryan, Jim; Leung, Naichung Conan
{\it The enumerative geometry of $K3$ surfaces and modular forms.}
J. Amer. Math. Soc. 13 (2000), no. 2, 371--410 (electronic). 


\bibitem{FV1}Friedl, Stefan; Vidussi, Stefano. {\it Twisted Alexander polynomials detect fibered 3-manifolds 
.}  arXiv:0805.1234 .

\bibitem{FV2}Friedl, Stefan; Vidussi, Stefano. {\it  Symplectic 4-manifolds with a free circle action.}   arXiv:0801.1513.

\bibitem{FM}Friedman, Robert ;  Morgan, John W.  {\it Algebraic surfaces and Seiberg-Witten invariants.}
 J. Algebraic Geom.  6  (1997),  no. 3, 445--479.

\bibitem{SF}Sakamoto, Koichi; Fukuhara, Shinji.{\it Classification of $T\sp{2}$-bundles over $T\sp{2}$.}
Tokyo J. Math. 6 (1983), no. 2, 311--327.

\bibitem{Ge} Geiges, Hansj\"org . {\it Symplectic structures on $T\sp 2$-bundles over $T\sp 2$.}
 Duke Math. J.  67  (1992),  no. 3, 539--555.

\bibitem {Go2} Gompf, Robert E.  {\it A new construction of symplectic manifolds.}
 Ann. of Math. (2)  142  (1995),  no. 3, 527--595.

\bibitem {Go5} Gompf, Robert E.  {\it Locally holomorphic maps yield symplectic structures.}
 Comm. Anal. Geom.  13  (2005),  no. 3, 511--525.
\bibitem{GS} Gompf, Robert E.; Stipsicz, Andr\'as I.
{\it $4$-manifolds and Kirby calculus. }
Graduate Studies in Mathematics, 20. American Mathematical Society, Providence, RI, 1999. xvi+558 pp.

\bibitem{IP5}  Ionel, Eleny-Nicoleta ;  Parker, Thomas H.  {\it The symplectic sum formula for Gromov-Witten invariants.}
 Ann. of Math. (2)  159  (2004),  no. 3, 935--1025.

\bibitem{IP4}  Ionel, Eleny-Nicoleta ;  Parker, Thomas H.  {\it Relative Gromov-Witten invariants.}
 Ann. of Math. (2)  157  (2003),  no. 1, 45--96.

\bibitem{J} Jabuka, Stanislav. {\it Symplectic surfaces and generic $J$-holomorphic structures on 4-manifolds. (English summary)}
Illinois J. Math. 48 (2004), no. 2, 675--685.

\bibitem{kod}Kodaira, Kunihiko.  {\it On the structure of compact complex analytic surfaces. I.}
 Amer. J. Math.  86  (1964), 751--798.

\bibitem{Le}  Lerman, Eugene . {\it Symplectic cuts.}
 Math. Res. Lett.  2  (1995),  no. 3, 247--258.


\bibitem{TJBL}Li, Bang-He; Li, Tian-Jun.{\it Symplectic genus, minimal genus and diffeomorphisms.} Asian J. Math. 6 (2002), 123--144.

\bibitem{TJL1}Li, Tian-Jun. {\it The Space of Symplectic Structures on closed 4-Manifolds.} AMS/IP Studies in Advanced Mathematics, V. 42  (2008),  259--273. (arxiv 0805.2931)

\bibitem{TJL2}Li, Tian-Jun. {\it Smoothly embedded spheres in symplectic {$4$}-manifolds.} Proc. Amer. Math. Soc. 2(1999), 609-613.

\bibitem {TJLL} Li, Tian-Jun ;  Liu, Ai-Ko. {\it  Uniqueness of symplectic canonical class, surface cone and symplectic
 cone of 4-manifolds with $b\sp +=1$.}
 J. Differential Geom.  58  (2001),  no. 2, 331--370.

\bibitem{LU}Li, Tian-Jun ;  Usher, Michael. {\it Symplectic forms and surfaces of negative square.}
 J. Symplectic Geom.  4  (2006),  no. 1, 71--91.


\bibitem {MW}McCarthy, John D. ;  Wolfson, Jon G.  {\it Symplectic normal connect sum.}
 Topology  33  (1994),  no. 4, 729--764.

\bibitem {M1}McDuff, Dusa. {\it From symplectic deformation to isotopy.}
 Topics in symplectic $4$-manifolds (Irvine, CA, 1996), 
 85--99, First Int. Press Lect. Ser., I, Int. Press, Cambridge, MA,  1998.

\bibitem {M2} McDuff, Dusa. {\it Notes on ruled symplectic $4$-manifolds.}
 Trans. Amer. Math. Soc.  345  (1994),  no. 2, 623--639.
                

\bibitem {M3} McDuff, Dusa. {\it Lectures on Gromov invariants for symplectic $4$-manifolds.}
With notes by Wladyslav Lorek.
NATO Adv. Sci. Inst. Ser. C Math. Phys. Sci., 488,  Gauge theory and symplectic geometry (Montreal, PQ, 1995), 
 175--210, Kluwer Acad. Publ., Dordrecht,  1997.

\bibitem{M4} McDuff, Dusa. {\it The local behaviour of holomorphic curves in almost complex
 $4$-manifolds.}
 J. Differential Geom.  34  (1991),  no. 1, 143--164.

\bibitem{M5} McDuff, Dusa. {\it Singularities and positivity of intersections of $J$-holomorphic
 curves.}
With an appendix by Gang Liu.
Progr. Math., 117,  Holomorphic curves in symplectic geometry, 
 191--215, Birkh\"auser, Basel,  1994.


\bibitem{McS}  McDuff, Dusa ;  Salamon, Dietmar. {\it $J$-holomorphic curves and symplectic topology.}
American Mathematical Society Colloquium Publications, 52. American Mathematical Society, Providence, RI,  2004. xii+669 pp. ISBN: 0-8218-3485-1 

\bibitem{PS}Pjatecki\u\i -\v Sapiro, I. I. ;  \v Safarevi\v c, I. R.  {\it Torelli's theorem for algebraic surfaces of type ${\rm K}3$.}
(Russian)  Izv. Akad. Nauk SSSR Ser. Mat.  35  (1971), 530--572. 


\bibitem{T}  Taubes, Clifford Henry. {\it Seiberg Witten and Gromov invariants for symplectic $4$-manifolds.}
Edited by Richard Wentworth.
First International Press Lecture Series, 2. International Press, Somerville, MA, 2000. vi+401 pp. ISBN: 1-57146-061-6


\bibitem{Ue}Ue, Masaaki. {\it Geometric $4$-manifolds in the sense of Thurston and Seifert
 $4$-manifolds. I.}
 J. Math. Soc. Japan  42  (1990),  no. 3, 511--540.

\bibitem{U} Usher, Michael . {\it The Gromov invariant and the Donaldson-Smith standard surface
 count.}
 Geom. Topol.  8  (2004), 565--610 (electronic).

\bibitem {Witten}  Witten, Edward . {\it Monopoles and four-manifolds.}
 Math. Res. Lett.  1  (1994),  no. 6, 769--796.
\end{thebibliography}
\end{document}